\documentclass[12pt]{article}
\usepackage{a4wide}
\usepackage{amssymb}
\usepackage{amsmath}
\usepackage[russian]{babel}

  \chardef\No=242
\newtheorem{theorem}{Theorem}{}
{}
{}
{}
{}

\begin{document}

\begin{center}
{\bf Corrections to the paper "The boundedness of certain sublinear operator in the weighted variable Lebesgue spaces".}\\
Rovshan A. Bandaliev \\
Institute of Mathematics and Mechanics of National Academy of Sciences of Azerbaijan\\
bandaliyev.rovshan@math.ab.az
\end{center}

In this paper author was proved the boundedness of the multidimensional Hardy type operator in weighted
Lebesgue spaces with a variable exponent. As an application we prove the boundedness of certain
sublinear operators on the weighted variable Lebesgue space. Note that the proof of multidimensional
Hardy type operator in weighted Lebesgue spaces with a variable exponent don't contained any mistakes.
But at the proving of the boundedness of certain sublinear operators on the weighted variable Lebesgue
space Georgian colleagues discovered a small but significant error in my paper, which was published in {\it Czechoslovak Math. Journal, 2010, 60 (135) (2010), 327--337.}  

This result is assigned as Theorem 5 in noted paper.
In other words, sufficient conditions for general weights ensuring the validity of the two-weight
strong type inequalities for some sublinear operator was found. In this theorem the inequality (9) isn't true.
In this note we give the details of the correct argument. We presume that the reader is familiar with the contents
and notation of our original paper. At the heart of our correction is the following Theorem which replaces Theorem 5.
The numeration in this note is remains as in [1].
\begin{theorem}
Let $1< \underline p\le p(x)\le \overline p< \infty$ for $x\in R^n$ and $\int_{R^n} \delta^{\frac {\underline p\,p(x)}{p(x)- \underline p}}\, dx< \infty,$ where $\delta\in [0, 1),$ and let $T$ be a sublinear operator acting boundedly from $L_{\underline p}(R^n)$ to $L_{p(x)}(R^n)$ such that, for any $f\in L_1(R^n)$ with compact support and $x\notin supp\; f$
$$
|Tf(x)|\le C\,\int\limits_{R^n} \frac{|f(y)|}{|x- y|^n}\,dy,
\leqno(3)
$$
where $C> 0$ is independent of $f$ and $x.$

Moreover, let $v(x)$ and $w(x)$ are weight functions on $R^n$ satisfying the conditions
$$
A_1= \sup\limits_{t> 0}\left(\int\limits_{|y|< t} [v(y)]^{-
\overline p^{\,'}}\;dy\right)^{\frac{\alpha}{\overline p{\,'}}}
\left\|\frac{w(\cdot)}{|\cdot|^n}\left(\int\limits_{|y|< |\cdot|}
[v(y)]^{- \overline p^{\,'}}\;dy\right)^{\frac{1- \alpha} {\overline
p^{\,'}}}\right\|_{L_{p(\cdot)}(|\cdot|> t)}< \infty, \leqno(4)
$$
$$
B= \sup\limits_{t> 0}\left(\int\limits_{|y|> t} [v(y)|y|^n]^{-
\overline p^{\,'}}\;dy\right)^{\frac{\beta}{\overline p{\,'}}}
\left\|w(\cdot)\left(\int\limits_{|y|> |\cdot|} [v(y)|y|^n]^{-
\overline p^{\,'}}\;dy\right)^{\frac{1- \beta} {\overline
p^{\,'}}}\right\|_{L_{p(\cdot)}(|\cdot|< t)}< \infty, \leqno(5)
$$
where $0< \alpha, \beta< 1.$

There exists $M> 0$ such that
$$
\sup\limits_{|x|/4< |y|\le 4\,|x|} w(y)\le M\, ess \,\inf_{x\in R^n}
v(x). \leqno(6)
$$

Then there exists a positive constant $C,$ independent of $f,$ such that for all $f\in L_{p(x), v}(R^n)$
$$
\|Tf\|_{L_{p(x), w}(R^n)}\le C \|f\|_{L_{p(x), v}(R^n)}.
$$
\end{theorem}

{\bf Proof.} For $k\in Z$ we define $E_k= \left\{x\in R^n:\; 2^k< |x|\le 2^{k+ 1}\right\},$ $E_{k,1}= \left\{x\in R^n:\; |x|\le 2^{k- 1}\right\},$ $E_{k,2}= \left\{x\in R^n:\;2^{k- 1}< |x|\le 2^{k+ 2}\right\}$ and $E_{k, 3}= \left\{x\in R^n:\; |x|> 2^{k+ 2}\right\}.$ Then $E_{k, 2}= E_{k- 1}\cup E_k\cup E_{k+ 1}$ and the multiplicity of the covering $E_{k,2}$ is equal to 3.

Given $f\in L_{p(x), v}(R^n),$ we have
$$
|Tf(x)|= \sum\limits_{k\in Z} |Tf(x)|\,\chi_{E_k}(x)\le
\sum\limits_{k\in Z} \left|Tf_{k,1}(x)\right|\,\chi_{E_k}(x)+
\sum\limits_{k\in Z} \left|Tf_{k,2}(x)\right|\,\chi_{E_k}(x)+
$$
$$
+ \sum\limits_{k\in Z} \left|Tf_{k,3}(x)\right|\,\chi_{E_k}(x) = T_1
f(x)+ T_2 f(x)+ T_3 f(x),
$$
where $\chi_{E_k}$ is the characteristic function of the set $E_k,$ $f_{k,i}= f\chi_{E_{k, i}},$ $i= 1, 2, 3.$

The estimates $\left\|T_1 f\right\|_{L_{p(x), w}(R^n)}\le C \|f\|_{L_{p(x), v}(R^n)}$ and $\left\|T_3 f\right\|_{L_{p(x), w}(R^n)}\le C \|f\|_{L_{p(x), v}(R^n)}$ is precisely and remains as in Theorem 5 in [1].

The mistake of author was in assuming that the inequality
$$
\left\|T_2 f\right\|_{L_{p(x), w}(R^n)}\le \|T\|_{L_{p(\cdot)}\left(R^n\right)}\, M\, \sum\limits_{k\in Z} \|f\,v\|_{L_{p(x)}\left(E_{k, 2}\right)}\le C_3\, \|f\,v\|_{L_{p(x)}\left(R^n\right)} \leqno(9)
$$
is holds. But after publication is discovered that the inequality (9) isn't true.

Now we reduce the correct variant of this inequality.
Since the operator $T$ is sublinear it suffices to prove that from $\|f\|_{L_{\underline p,\,v}\left(R^n\right)}\le 1$ implies $\displaystyle{\int\limits_{R^n} w^{p(x)}(x)\left[\sum\limits_{k\in Z} \left|Tf_{k, 2}\right|\,\chi_{E_k} \right]^{p(x)}\,dx\le C},$ where $C> 0$ is independent on $k\in Z$ (see [2]).

We have
$$
\int\limits_{R^n} w^{p(x)}(x)\left[\sum\limits_{k\in Z} \left|Tf_{k, 2}\right|\chi_{E_k} \right]^{p(x)}\,dx=
\sum\limits_{k\in Z}\, \int\limits_{E_k} \, \left(\left|Tf_{k, 2}\right|\,w(x)\right)^{p(x)}\,dx.
$$
By virtue of boundedness of operator $T$ and condition (6), we have
$$
\sum\limits_{k\in Z}\, \int\limits_{E_k} \, \left(\left|Tf_{k, 2}\right|\, w(x)\right)^{p(x)}\,dx=
\sum\limits_{k\in Z}\, \int\limits_{E_k} \, \left(\frac{\left|Tf_{k, 2}\right|} {C_1\,\|f_{k,2}\|_{L_{\underline p} \left(R^n\right)}}\right)^{p(x)}\,\left(C_1\,\|f_{k, 2} \|_{L_{\underline p}}\, w(x)\right)^{p(x)}\,dx
$$
$$
\le C_2\,\sum\limits_{k\in Z}\,\sup\limits_{x\in E_k}\left(\|f_{k, 2} \|_{L_{\underline p}\left(R^n\right)}\,
w(x)\right)^{p(x)}\,\int\limits_{R^n} \, \left(\frac{\left|Tf_{k, 2}\right|} {C_1\,\|f_{k,2}\|_{L_{\underline p} \left(R^n\right)}}\right)^{p(x)}\,dx
$$
$$
\le C_2\,\sum\limits_{k\in Z}\,\sup\limits_{x\in E_k}\left(\|f_{k, 2} \|_{L_{\underline p} \left(R^n\right)}\,w(x)\right)^{p(x)}= C_2\,\sum\limits_{k\in Z}\,\sup\limits_{x\in E_k} \left(\|f\|_{L_{\underline p}\left(E_{k,2}\right)}\,w(x)\right)^{p(x)}
$$
$$
= C_2\,\sum\limits_{k\in Z}\,\sup\limits_{x\in E_k} \left(\|f\,w(x)\|_{L_{\underline p}\left(E_{k,2}\right)} \right)^{p(x)}\le C_3 \sum\limits_{k\in Z}\,\sup\limits_{x\in E_k} \left(\|f\,\inf\limits_{x\in E_{k,2}}\,v(x)\|_{L_{\underline p}\left(E_{k,2}\right)}\right)^{p(x)}
$$
$$
\le C_3 \sum\limits_{k\in Z}\,\sup\limits_{x\in E_k} \left(\|f\,v\|_{L_{\underline p} \left(E_{k,2}\right)} \right)^{p(x)}= C_3 \sum\limits_{k\in Z}\,\sup\limits_{x\in E_k} \left(\|f\|_{L_{\underline p,\,v} \left(E_{k,2}\right)}\right)^{p(x)}
$$
$$
= C_3 \sum\limits_{k\in Z} \left(\|f\|_{L_{\underline p,\,v}\left(E_{k,2}\right)}\right)^{\inf\limits_{x\in E_k} p(x)}\le C_3 \sum\limits_{k\in Z} \left(\|f\|_{L_{\underline p,\,v}\left(E_{k,2}\right)}\right)^{\underline p}.
$$
Further, we obtain
$$
\sum\limits_{k\in Z} \left(\|f\|_{L_{\underline p,\,v}\left(E_{k,2}\right)}\right)^{\underline p}=
\sum\limits_{k\in Z} \left(\int\limits_{E_{k, 2}} |f(x)\,v(x)|^{\underline p}\,dx\right)=
$$
$$
\sum\limits_{k\in Z} \left(\int\limits_{E_{k- 1}}+ \int\limits_{E_k}+ \int\limits_{E_{k+ 1}}\right) |f(x)\,v(x)|^{\underline p}\,dx=
$$
$$
3\,\left(\int\limits_{R^n} |f(x)\,v(x)|^{\underline p}\,dx\right)= 3\; \|f\|_{L_{\underline p,\,v}\left(R^n\right)}^{\underline p}\le 3.
$$
Thus $\left\|Tf_{k, 2}\right\|_{L_{p(\cdot),\,w}\left(R^n\right)}\le C\, \|f\|_{L_{\underline p,\,v}\left(R^n\right)}$ and take into account the condition $\displaystyle{\int\limits_{R^n} \delta^{\frac {\underline p\,p(x)}{p(x)- \underline p}}\, dx< \infty},$ we obtain
$$
\left\|Tf_{k, 2}\right\|_{L_{p(\cdot),\,w}\left(R^n\right)}\le C\, \|f\|_{L_{\underline p,\,v}\left(R^n\right)}\le
C \|f\|_{L_{p(\cdot),\,v}\left(R^n\right)}.
$$
The proof of Theorem 5 is complete.

{\bf Acknowledgement.} The authors thank Prof. Dr. V. Kokilashvili who called my attention to the gap in the
proof.
\vspace{3mm}

\begin{center}
{\bf References}
\end{center}

[1] R. A. Bandaliev, {\it The boundedness of certain sublinear operator in the weighted variable Lebesgue spaces},
Czechoslovak Math. J. {\bf 60}(2010), no. 2, 327-337.

[2] L. Diening, P. Harjulehto, P. H\"{a}st\"{o} and M. R\.{u}\v{z}i\v{c}ka,
{\it Lebesgue and Sobolev spaces with variable exponents,} Springer Lecture Notes,
v.2017, Springer-Verlag, Berlin, 2011.

\end{document}